\documentclass[draft,12pt, reqno]{amsart}
\usepackage{amssymb, amsmath}
\usepackage{url}
\usepackage{enumerate}
\usepackage{mathabx}

\newtheorem{theorem}{Theorem}[section]
\newtheorem{corollary}[theorem]{Corollary}
\newtheorem{prop}[theorem]{Proposition}
\newtheorem{lemma}[theorem]{Lemma}

\theoremstyle{remark}

\theoremstyle{definition}

\newtheorem{example}[theorem]{Example}

\numberwithin{equation}{section}
\numberwithin{theorem}{section}

\newcommand{\sgn}{{\rm sgn}}

\newcommand{\bv}{{\mathcal BV}}
\newcommand{\bvRbar}{{\mathcal BV}([-\infty,\infty])}

\newcommand{\N}{{\mathbb N}}
\newcommand{\R}{{\mathbb R}}
\newcommand{\C}{{\mathbb C}}

\newcommand{\fn}{\!:\!}

\newcommand{\fhat}{\hat{f}}

\providecommand{\abs}[1]{\lvert#1\rvert}
\providecommand{\norm}[1]{\lVert#1\rVert}

\newcommand{\Lone}{L^1(\R)}
\newcommand{\intab}{\int_a^b}

\newcommand{\intinf}{\int_{-\infty}^{\infty}}

\newcommand{\sumN}{\sum_{n=1}^N}
\newcommand{\reg}{{\mathcal R}}

\begin{document}
\date{Preprint March 25, 2022.  To appear in {\it Mathematica Slovaca}.}
\title[Fourier transform inversion]
{Fourier transform inversion: Bounded variation, polynomial growth, Henstock--Stieltjes integration}
\author{Erik Talvila}
\address{Department of Mathematics \& Statistics\\
University of the Fraser Valley\\
Abbotsford, BC Canada V2S 7M8}
\email{Erik.Talvila@ufv.ca}
\thanks{An anonymous referee provided helpful comments and several important references.}

\begin{abstract}
In this paper we prove pointwise and distributional Fourier transform inversion theorems for functions
on the real line
that are locally of bounded variation, while in a neighbourhood of infinity are
Lebesgue integrable or have polynomial growth.  We also allow the Fourier transform
to exist in the principal value sense.  A function is called regulated
if it has a left limit and a right limit at each point.
The main inversion theorem is obtained by
solving the differential equation $df(t)-i\omega f(t)=g(t)$ for a regulated function $f$,
where $\omega$ is a complex number with positive imaginary part.
This is done using the
Henstock--Stieltjes integral.  This is an integral defined with Riemann sums
and a gauge. Some variants of
the integration by parts formula are also proved for this integral.
When the function is of polynomial growth its Fourier transform exists in a distributional
sense, although the inversion formula only involves integration of functions and returns
pointwise values.
\end{abstract}

\maketitle

\section{Introduction}\label{sectionintroduction}

If $f\fn\R\to\R$ then its Fourier transform is $\fhat(s)=\intinf e^{-ist}f(t)\,dt$.  
A sufficient condition for existence
of $\fhat$ on $\R$ is that $f\in \Lone$;  and then $\fhat$ is uniformly continuous
on $\R$.  Under the same hypothesis, the Riemann--Lebesgue lemma
says $\fhat(s)$  has limit $0$ as $\abs{s}\to\infty$.  We also use the fact that if $f$ is positive
and decreases to $0$ then the integral $\int_0^\infty e^{-ist}f(t)\,dt$ exists for $s\not=0$.
However, under this condition we need not have $\fhat(s)\to 0$ as $s\to\infty$.

The process of recovering $f$ from $\fhat$ is known as Fourier inversion.
A basic theorem is that of Jordan.
If $f\in\Lone$ and is also of 
bounded variation on the compact interval $[a,b]$ then for each $x\in(a,b)$
\begin{equation}
\frac{f(x-)+f(x+)}{2}=\frac{1}{2\pi}\lim_{T\to\infty}\int_{-T}^T e^{ixs}\fhat(s)\,ds.\label{inversionformula}
\end{equation}
Here, the inversion is given by a principal value integral.  The usual proofs use  the Riemann--Lebesgue
lemma and properties of the
Dirichlet kernel.  See,
for example, \cite{bachman},
\cite{benedetto} and \cite[\S1.9]{titchmarsh}.
These references also supply a suitable background on Fourier transforms.

We prove Fourier inversion theorems under the following assumptions on $f$.  If the
inversion is to hold at $x\in(a,b)$ then $f$ is assumed to be of bounded variation on $[a,b]$.
On each compact interval that does not intersect $(a,b)$, the function $f$ is locally Lebesgue
integrable or has a principal value integral.  In a neighbourhood of infinity, $f$ is Lebesgue
integrable, is of bounded variation, or is asymptotic to a polynomial.  In the latter two cases,
the Fourier transform may only exist as a distribution and yet the inversion only involves
integration of functions and yields pointwise values.

There are of course many criteria besides local bounded variation in the literature applied 
to behaviour of $f$ at the inversion
point.  For example, Dini's test or degrees of differentiability. See \cite{bachman} and \cite{titchmarsh}.
Koekoek \cite{koekoek} has provided a simple proof of the
inversion theorem that applies at points of differentiability or for Lipschitz or H\"older 
continuous functions.

The outline of the paper is as follows.
For a regulated function (one having a left limit and a right limit
at each point) the following identity is proved in Lemma~\ref{lemmahsrep}:
\begin{equation}
\left[\frac{f(x-)+f(x+)}{2}\right] = 
\intinf H(x-t) e^{i\omega(x- t)}[df(t)-i\omega f(t)\,dt].\label{regulatedidentity}
\end{equation}
Here $H$ is the Heaviside step function.
The above identity is proved using the Henstock--Stieltjes integral, which is described in 
Section~\ref{sectionHSintegral}.
Some formulas related to integration by parts are proved here.
An integral representation of the Heaviside step function is
given in Lemma~\ref{lemmaperron1}.  Inserting this in \eqref{regulatedidentity} and changing the
orders of integration
completes the proof of Jordan's theorem
(Theorem~\ref{theoremjordan}). Pringsheim's theorem, for which $f$ is of bounded variation with limit
$0$ at $\pm\infty$ but need not be in $\Lone$,
follows as Corollary~\ref{corollarypringsheim}.
In Pringsheim's theorem, when the
condition that $f$ vanishes at infinity is dropped, the Fourier transform may only exist as a 
distribution (generalised function).  However, a pointwise inversion theorem, in which only functions need
to be
integrated, is given in Theorem~\ref{theorembvlimit}.  Similarly in Theorem~\ref{theorempolynomial}
 for functions  that are asymptotic to a
polynomial in a
neighbourhood of infinity.  Functions for which the Fourier transform exists in a 
principal value sense are considered in Theorem~\ref{theorempv}.

The Henstock--Stieltjes integral is defined using a gauge that controls function
evaluations in Riemann sums.  Unlike the Riemann--Stieltjes integral, integration
is defined over the entire real line without resorting to improper integrals.  Also,
the function and integrator can have coincident jump discontinuities.  In Example~\ref{exampleRS}
it is shown how to compute
\eqref{regulatedidentity} as a Henstock--Stieltjes integral and that the corresponding
Riemann--Stieltjes integral can fail to exist,
even allowing for improper integrals.

Our results assume the function is of bounded variation in a neighbourhood of the point $x$ in 
\eqref{inversionformula}.  For other Fourier analysis results for functions of bounded
variation, see \cite{liflyand2018}, \cite{liflyand2019} and \cite{moricz2004}.  This latter paper uses improper
Riemann--Stieltjes integrals.

Some of our results can be proved using summability techniques for distributions
developed by R.~Estrada and J.~Vindas.
See \cite{pilipovicstankovicvindas}, Section~II.5.2 and references therein.

\section{Preliminary results}\label{sectionpreliminary}

The two main parts of our proof of the inversion theorem are a contour integral for the Heaviside step
function and a representation of a function as a Henstock--Stieltjes integral.

The Heaviside step function is
\begin{equation}
H(x)=\left\{\begin{array}{cl}
1, \text{ if } x>0\\
1/2, \text{ if } x=0\\
0, \text{ if } x<0.
\end{array}
\right.
\end{equation}
And, the signum function is ${\rm sgn}(x)=2H(x)-1$.
\begin{lemma}\label{lemmaperron1}
Let $p$ be a real number and let $\omega$ be a complex number with positive imaginary part.  Then
\begin{equation}
\frac{1}{2\pi i}\intinf \frac{e^{ips}}{s-\omega}\,ds =H(p)e^{ip\omega}.\label{Hrep}
\end{equation}
Carrying out the same calculation with $\omega=0$ shows that $\int_0^\infty [\sin(px)/x]\,dx=(\pi/2)
{\rm sgn}(p)$.
\end{lemma}
If $p\not=0$,
the integral in \eqref{Hrep} converges conditionally.
If $p=0$ or $\omega$ is real, the integral exists in the principal value sense.
It is sometimes called Perron's lemma, although it seems to have been first derived by
Cauchy (cf. \cite[p.~123]{whittaker}).
We quote some estimates from the proof in \cite{talvilafouriermaa} in order to justify
some limit operations in the proof of Theorem~\ref{theoremjordan}.

Suppose $g\in L^1(\R)$ and $\omega\in\C$.  If
$f$ is absolutely continuous and is a solution of the differential equation
$f'(t)-i\omega f(t)=g(t)$ then 
\begin{align}
f(x)  &=  f(0)e^{i\omega x}+e^{i\omega x}\int_{0}^x e^{-i\omega t}g(t)\,dt\notag\\
  &=  f(0)e^{i\omega x}+ \int_0^\infty H(x-t)e^{i\omega(x-t)}[f'(t)-i\omega f(t)]\,dt.\label{ACDE}
\end{align}
Representations of this type were used in \cite{talvilafouriermaa} to prove a Fourier inversion theorem.
See Corollary~\ref{corollaryac} below.

If $-\infty\leq a<b\leq\infty$ then
the function $f\fn(a,b)\to\R$ is regulated if for each $x\in(a,b]$ and 
each $y\in[a,b)$ the one-sided limits exist;
$\lim_{t\to x^-}f(t)=f(x-)\in\R$ and $\lim_{t\to y^+}f(t)=f(y+)\in\R$.  
If $f$ is defined on $(a,b)$ we will
assume its domain has been extended to $[a,b]$ with the above limits.  Regulated functions are bounded
and have at most a countable number of finite jump discontinuities.  They are thus measurable and locally
summable.  Denote the space of regulated functions by $\reg([a,b])$.
For more
on regulated functions see \cite{mcleod} or \cite{monteiroslaviktvrdy}.

An analogue of \eqref{ACDE} for regulated functions is the following lemma.

\begin{lemma}\label{lemmahsrep}
Let $-\infty<a<x<b<\infty$. (a)  Let $\omega$ be a complex number with positive imaginary part.
Let $f\in L^1((-\infty,b])$ and regulated on $[-\infty,b]$.  Then
\begin{equation}
\frac{f(x-)+f(x+)}{2} = 
\intinf H(x-t) e^{i\omega(x- t)}[df(t)-i\omega f(t)\,dt].\label{hsrep}
\end{equation}
(b) Let $\omega$ be a complex number.  Let $f$ be regulated on $[a,b]$.  Then
\begin{equation}
\frac{f(x-)+f(x+)}{2}  =  \intab H(x-t)e^{i\omega(x-t)}[df(t)-i\omega f(t)\,dt]
 +f(a)e^{i\omega(x-a)}.\label{hsrepfinite}
\end{equation}
\end{lemma}

The proof uses the Henstock--Stieltjes integral and variants on the integration by parts formula.
See section~\ref{sectionHSintegral} and \ref{lemmaproofhsrep}.  The lemma is
false for Riemann--Stieltjes integrals.  See Example~\ref{exampleRS}.

For each $N\in\N$ a partition of $[a,b]$ is a collection of points $a=x_0<x_1<\ldots<x_N=b$.
If $-\infty\leq a<b\leq\infty$ then
function $f$ is of bounded variation on $[a,b]$ if the supremum of $\sumN\abs{f(x_n)-f(x_{n-1})}$ is bounded, 
where the supremum is taken over all partitions of $[a,b]$.
In this case the supremum over all partitions is labeled $V\!f$.  
Functions of bounded variation are regulated.  
For more on functions of bounded variation see \cite{appell} or \cite{mcleod}.

\section{Inversion theorem}\label{sectionft}

We now use formula \eqref{hsrep} and the integral representation of the Heaviside step function
in Lemma~\ref{lemmaperron1} to prove an inversion formula for functions of bounded variation.
Notice that the proof does not use the Riemann--Lebesgue lemma or the Dirichlet kernel.

\begin{theorem}[Jordan]\label{theoremjordan}
Let $f\in\bvRbar\cap L^1(\R)$.  Then for each $x\in\R$,
\begin{equation}
\frac{f(x-)+f(x+)}{2}=\frac{1}{2\pi}\lim_{T\to\infty}\int_{-T}^T e^{ixs}\fhat(s)\,ds.\label{inversion}
\end{equation}
\end{theorem}

\begin{proof}
With Lemma~\ref{lemmaperron1} and Lemma~\ref{lemmahsrep}  we have
\begin{eqnarray}
\frac{f(x-)+f(x+)}{2} & = & \intinf H(x-t)e^{i\omega(x-t)}[df(t)-i\omega f(t)\,dt]\label{repformula}\\
 & = & \frac{1}{2\pi i}\intinf\lim_{T\to\infty}\int_{-T}^T \frac{e^{i(x-t)s}}{s-\omega}\,ds[df(t)-i\omega f(t)\,dt]\notag\\
 & = & \frac{1}{2\pi i}\lim_{T\to\infty}\int_{-T}^T\frac{e^{ixs}}{s-\omega}\intinf 
e^{-ist}[df(t)-i\omega f(t)\,dt]\,ds.\notag
\end{eqnarray}
Since $f\in\Lone $ and is of bounded variation,
dominated convergence and then the Fubini--Tonelli theorem allow us to bring the limit outside the
integral and then change the orders of integration.  See equation (14) in \cite{talvilafouriermaa}.  
When $t\leq x$ this shows
$\abs{I_T}\leq 2\pi$.  The magnitude of the residue is 
$e^{-(x-t)\eta}\leq 1$, where $\eta$ is the imaginary part of $\omega$.
When $t>x$ we get the same estimate on $I_T$ when the half circle below the real axis is used.
Now the residue is zero.
Hence, we can bring the $T$ limit outside the integral.
And, $\abs{\frac{e^{i(x-t)s}}{s-\omega}}\leq 1/\eta$.
We can then change the orders of integration.

Now consider $\intinf e^{-ist}df(t)$.  Since $f\in\Lone $ and is of bounded variation it
necessarily has limit $0$ at $\pm\infty$.  Dirichlet's test then shows existence of the integral
for each $s\in\R$.
Defining this as a Henstock--Stieltjes
integral requires giving the function $t\mapsto e^{-ist}$ values at $t=\pm\infty$.
It is immaterial how this is done because every $\gamma$-fine tagged partition
of $[-\infty,\infty]$ will contain a term $e^{-is\infty}[f(\infty)-f(x_{N-1})]$.  No
matter what value is given to $e^{-is\infty}$ this term can be made as small as desired by
taking $x_N$ large enough.  Similarly at $-\infty$.  This lets us write
\begin{align*}
&\intinf e^{-ist}df(t)  =  \lim_{\substack{Y\to\infty\\X\to-\infty}}\int_X^Y e^{-ist}df(t)\\
 & =  \lim_{Y\to\infty}\left[e^{-isY}f(Y)\right] - \lim_{X\to-\infty}\left[e^{-isX}f(X)\right] 
-\lim_{\substack{Y\to\infty\\X\to-\infty}}\int_X^Y f(t)\,d\!\left[e^{-ist}\right]\\
 & =  is\intinf e^{-ist}f(t)\,dt,
\end{align*}
where
the integration by parts formula \eqref{parts} and Proposition~\ref{propB'} are applied on
each finite interval 
$[X,Y]$.
\end{proof}
Since the exponential function is continuous the last integral can also be evaluated
as the limit of a Riemann--Stieltjes integral over the compact interval $[X,Y]$ and
then limits taken.

An anonymous referee has pointed out that:
``Theorem~\ref{theoremjordan} is essentially contained in Theorem~15 from
\cite{estradavindas2007} when applied with any $1 < p < \infty$, 
since under the author hypothesis the Fourier transform under consideration is $O(1/|s|)$
 at infinity. Alternatively, Theorem~\ref{theoremjordan} follows at once by 
combining Theorem~13 in \cite{estradavindas2007} with the distributional version of 
Littlewood's Tauberian theorem from Theorem~4.1 in \cite{estradavindas2013}."

We now present three corollaries that follow with different assumptions on $f$.  In 
Corollary~\ref{corollaryac}
we note that if $f$ is absolutely continuous and $f'\in\Lone$ then $f\in\bvRbar$.
Usual methods with the Dirichlet kernel (for example, \cite[\S1.9]{titchmarsh}) in 
Corollary~\ref{corollary1} show
that bounded variation is only a local requirement.  Corollary~\ref{corollarypringsheim}
is similar to results due to Pringsheim.  See
\cite[p.~15]{titchmarsh}.

\begin{corollary}\label{corollaryac}
Let $f\in\Lone$ such that $f$ is absolutely continuous and $f'\in\Lone$.  Then for each $x\in\R$
\begin{equation}
f(x)=\frac{1}{2\pi}\lim_{T\to\infty}\int_{-T}^T e^{ixs}\fhat(s)\,ds.\label{acexample}
\end{equation}
\end{corollary}
This was also proved in \cite{talvilafouriermaa}, where the assumptions allowed a proof with the
Lebesgue rather than Henstock--Stieltjes integral.

\begin{example}\label{examplelog}
Define $f$ by
$$
f(x)=\left\{\begin{array}{cl}
0, & x\leq 0\\
1/\log(x), & 0<x\leq 1/e\\
-1/(ex)^2, & x\geq 1/e.
\end{array}
\right.
$$
Then $f$ is absolutely continuous, the conditions of the corollary hold and we get \eqref{acexample}
at each point.  Note that $f$ is not differentiable at $0$ or $1/e$ and that $f$ is not Lipschitz
or H\"older continuous at the origin.  Hence, Koekoek's \cite{koekoek} result does not apply at the
origin.
\end{example}

\begin{corollary}\label{corollary1}
Let $f\in\bv([a,b])$ for some $-\infty<a<b<\infty$ and $f\in L^1((-\infty,a])$ and $f\in L^1([b,\infty))$. 
Let $x\in(a,b)$.  Then
\begin{equation}
\frac{f(x-)+f(x+)}{2}=\frac{1}{2\pi}\lim_{T\to\infty}\int_{-T}^T e^{ixs}\intab e^{-ist}f(t)\,dt\,ds.
\label{localinversion}
\end{equation}
And,
$$
\lim_{\substack{B\to\infty\\A\to-\infty}}\int_A^B e^{ixs}\intinf e^{-ist} 
f(t)\chi_{(-\infty,a]\cup[b,\infty)}(t)\,dt\,ds=0.
$$
\end{corollary}
\begin{proof}
Applying Theorem~\ref{theoremjordan} to $f\chi_{(a,b)}$ gives \eqref{localinversion}.

Let $A<0$ and $B>0$.
By the Fubini--Tonelli theorem
$$
\int_A^B e^{ixs}\int_b^\infty e^{-ist}f(t)\,dt\,ds=i\int_{b-x}^\infty \frac{f(t+x)e^{-iBt}}{t}\,dt
-i\int_{b-x}^\infty \frac{f(t+x)e^{-iAt}}{t}\,dt.
$$
Using the Riemann--Lebesgue lemma,
this tends to $0$ as $A\to-\infty$ and $B\to\infty$.  Similarly for integration over $(-\infty,a]$.
\end{proof}

\begin{example}\label{examplesingularfunction}
Choose
$f\fn[a,b]\to\R$ to be an increasing function whose derivative vanishes almost
everywhere.
For example, the Cantor--Lebesgue function.  Then $f$ is of bounded variation but is not absolutely 
continuous.  The formula in
Corollary~\ref{corollary1} gives $f(x)$ at each $x\in(a,b)$.  The results in \cite{koekoek} do not apply.
\end{example}

\begin{corollary}[Pringsheim]\label{corollarypringsheim}
Let $f\in\bvRbar$ such that $\lim_{\abs{x}\to\infty}f(x)=0$.  Then $\fhat(s)$ exists for $s\not=0$.  
And, for each $x\in\R$,
$$
\frac{f(x-)+f(x+)}{2}=\frac{1}{2\pi}\lim_{\substack{T\to\infty\\S\to 0^+}}\int_{S<\abs{s}<T} 
e^{ixs}\intinf e^{-ist}f(t)\,dt\,ds.
$$
Suppose the function $x\mapsto e^{-iax}f(x)$ is of bounded variation for some $a\in\R$.  Then, 
for each $x\in\R$,
$$
\frac{f(x-)+f(x+)}{2}=\frac{1}{2\pi}\lim_{\substack{T\to\infty\\S\to a^+}}\int_{S<\abs{s-a}<T} 
e^{ixs}\intinf e^{-ist}f(t)\,dt\,ds.
$$
\end{corollary}
For other proofs, see \cite{rieszlivingston} and \cite[\S1.10]{titchmarsh}.
See the latter 
for a reference to Pringsheim's original paper.
\begin{proof}
Existence of $\fhat$ is by Dirichlet's test.

Let $-\infty<a<x<b<\infty$.
As in Corollary~\ref{corollary1} we have \eqref{localinversion}.

Integrating by parts and using the Fubini--Tonelli theorem,
\begin{align*}
&\int_{S<\abs{s}<T} 
e^{ixs}\int_b^\infty e^{-ist}f(t)\,dt\,ds\\
&=\qquad \int_{S<\abs{s}<T} \frac{e^{ixs}}{is}\left\{e^{-isb}f(b)+\int_b^\infty
e^{-ist}\,df(t)\right\}\,ds\\
&=\qquad -2f(b)\int_{(b-x)S}^{(b-x)T}\sin(s)\,\frac{ds}{s}-2\int_b^\infty
\int_{(t-x)S}^{(t-x)T}\sin(s)\,\frac{ds}{s}\,df(t).
\end{align*}
The integral $\int_\alpha^\beta\sin(s)\,ds/s$ is bounded independent of its limits of integration.
It is evaluated in Lemma~\ref{lemmaperron1}.  Using dominated convergence we can then
take limits as $S\to 0^+$ and $T\to\infty$ and this gives
$$
-\pi\left\{f(b)+\int_b^\infty\,df(t)\right\}=0.
$$
The integral over $(-\infty,a)$ is handled similarly.
\end{proof}

\begin{example}\label{examplebv0}
(a) Let $f(x)={\rm sgn}(x)/\log\abs{x}$ for $\abs{x}>e$ and $f(x)=0$, otherwise.  For $s\not=0$, integrate
by parts:
\begin{eqnarray*}
\frac{i\fhat(s)}{2} & = & \int_e^\infty\frac{\sin(st)}{\log(t)}\,dt
= \frac{\cos(es)}{s}-\frac{1}{s}\int_e^\infty\frac{\cos(st)}{t\log^2(t)}\,dt\\
 & \sim &  \frac{\cos(es)}{s}
\text{ as } s\to\infty.
\end{eqnarray*}
The last line is by
Riemann--Lebesgue lemma.
Hence, $\fhat\not\in\Lone$.  The principal value limit $T\to\infty$ in Corollary~\ref{corollarypringsheim}
is required when $\abs{x}=e$ but for other values of $x$ the inversion integral exists as a Henstock--Kurzweil
integral.

Let $0<s<1/e$.
The second mean value theorem gives
$$
\left|\int_{1/s}^\infty\frac{\cos(st)}{t\log^2(t)}\,dt\right|=\frac{1}{\log^2(s)}\left|\int_1^{\xi}
\frac{\cos(t)}{t}\,dt\right|\leq \frac{C}{\log^2(s)}
$$
for some $\xi\geq 1$,
where $C=\sup_{x\geq 1}\abs{\int_1^x\cos(t)dt/t}$. 

By Taylor's theorem, there is $es\leq\eta(s,t)\leq 1$ such that
$$
\int_e^{1/s}\frac{\cos(st)}{t\log^2(t)}\,dt  =  \int_e^{1/s}\frac{dt}{t\log^2(t)}
-s\int_e^{1/s}\frac{\sin(\eta(s,t))}{t\log^2(t)}\,dt
$$
so that
$$
\left|\int_e^{1/s}\frac{\cos(st)}{t\log^2(t)}\,dt-1-\frac{1}{\log(s)}\right|\leq s,
$$
the last estimate being by dominated convergence.  These results show that
$i\fhat(s)/2\sim-1/(s\abs{\log\abs{s}})$ as $s\to0$.
Then $\fhat$ is not integrable in any right neighbourhood of the origin but 
$\sin(xs)\fhat(s)$ is bounded as $s\to 0^+$ so the principal value limit
$S\to 0^+$ in Corollary~\ref{corollarypringsheim} exists.

(b) An example of a function $f\in\bv([-\infty,\infty])$ that has limit $0$ at $\pm\infty$
but is not in $\Lone$ is $f(x)=\sin(x^{1/2})/x^{2/3}$ for $x>1$ and $f(x)=0$ for $x\leq 1$.
The integral of $f$ then converges conditionally.
Fourier transforms of such functions have been studied by J.H. Arredondo, F.J. Mendoza 
and A. Reyes
\cite{arredondo}.  Asymptotics of the Fourier--Laplace transform of such functions
is considered in \cite{brouckedebruyneindas}.
\end{example}

\section{Distributional transforms with pointwise inversion}

If $f$
is of bounded variation then it has a limit at $\pm\infty$.  Subtracting off this limit leaves a
function satisfying the conditions of Corollary~\ref{corollarypringsheim}.  The functions
$H(x)$ and $H(-x)$ have distributional Fourier transforms so $\fhat$ exists in this sense as well.
However, we can obtain a pointwise inversion
formula as with \eqref{inversionformula}.  This is done in the following theorem.
Theorem~\ref{theorempolynomial} considers the case when $f$ is asymptotic to a polynomial.

An anonymous referee has pointed out that:
``Theorem~\ref{theorembvlimit} and Theorem~\ref{theorempolynomial} also follow from Theorem~13 
and Theorem~5 in \cite{estradavindas2007} because trigonometric integrals (or more 
precisely principal value special evaluations in the sense of \cite{estradavindas2007}) 
associated to the second sum are easily seen to be Ces\`aro summable to 
$H(x)p_+(x) + H(-x)p_-(x)$ and moreover if $G$ is a primitive of $e^{ixs}{\hat g}(s)$, then
$$
\lim_{s\to 0^+}\left[G(s)-G(-s)\right]=2i\lim_{s\to0^+}\intinf\frac{\sin(st)}{st}\,df(t-x)
=2i\intinf df(t-x)=0."
$$
\begin{theorem}\label{theorembvlimit}
Let $f\in\bvRbar$, where we define $f(\infty)=\lim_{x\to\infty}f(x)$ and $f(-\infty)=\lim_{x\to-\infty}f(x)$.
Let $g(x)=f(x)-f(\infty)H(x)-f(-\infty)H(-x)$.
Then the conditions of Corollary~\ref{corollarypringsheim} apply to $g$.  Let $\delta$ be the Dirac
distribution.  Then
\begin{equation}
\fhat(s)={\hat g}(s)+f(\infty)\left[\pi\delta(s)+\frac{1}{i s}\right] +
f(-\infty)\left[\pi\delta(s)-\frac{1}{i s}\right].\label{fhatdelta}
\end{equation}
And, for each $x\in\R$,
\begin{equation}
\frac{f(x-)+f(x+)}{2}=\frac{1}{2\pi}\lim_{\substack{T\to\infty\\S\to 0^+}}\int_{S<\abs{s}<T} 
e^{ixs}\hat{g}(s)\,ds+f(\infty)H(x)+f(-\infty)H(-x).
\label{bvinversion}
\end{equation}
\end{theorem}
The distribution $T=1/s$ is given as $\langle T,\phi\rangle=\int_0^\infty\frac{\phi(x)-\phi(-x)}{x}\,dx$
for each test function $\phi$.  See \cite[p.~25]{gelfandshilov}. 

Although $\fhat$ is in general a distribution, notice that the inversion formula \eqref{bvinversion}
only requires integration of functions.

Limits of integrals,
such as in \eqref{bvinversion}, can be written in principal value
notation.  See Definition~5.8, page~277, in \cite{pilipovicstankovicvindas}.

\begin{proof}
The following distributional Fourier transforms are well-known,
\begin{equation}
{\hat H}(s)  =  \pi\delta(s)+\frac{1}{is},\qquad
{\widehat \sgn}(s)  =  \frac{2}{is}.\label{Hsignumft}
\end{equation}
See,  for example, \cite[p.~172, 174]{gelfandshilov}.  These give  the formula for $\fhat$ in the statement
of the theorem.

Notice that $g\in\bvRbar$ and $\lim_{\abs{x}\to\infty}g(x)=0$.  Applying Corollary~\ref{corollarypringsheim}
to $g$ and rearranging for $f$ gives the inversion formula \eqref{bvinversion}.
\end{proof}

\begin{example}
Let $f(x)=\arctan(x)$.  Then $g(x)=\arctan(x)-(\pi/2){\rm sgn}(x)$.  Now
integrate by parts to get
$$
{\hat g}(s)  =  -2i\int_0^\infty\left[\arctan(t)-\frac{\pi}{2}\right]\sin(st)\,dt
  =  \frac{2i}{s}\left\{\frac{\pi}{2}-\int_0^\infty\frac{\cos(st)}{1+t^2}\,dt\right\}.
$$
The cosine integral above can be evaluated with a contour integral as in the proof of
Lemma~\ref{lemmaperron1}.  Its value is $\pi e^{-\abs{s}}/2$ and then
${\hat g}(s)=i\pi(1-e^{-\abs{s}})/s$.  Equation~\eqref{fhatdelta} then gives
$\fhat(s)=-i\pi e^{-\abs{s}}/s$.
The inversion formula gives
\begin{eqnarray*}
\frac{1}{2\pi}\lim_{\substack{T\to\infty\\S\to 0^+}}\int_{S<\abs{s}<T} 
e^{ixs}\hat{g}(s)\,ds  & = &  \int_0^\infty e^{-s}\sin(xs)\,\frac{ds}{s}-\int_0^\infty \frac{
\sin(xs)}{s}\,ds\\
 & = & \arctan(x) -\frac{\pi}{2}{\rm sgn}(x).
\end{eqnarray*}
The principal value inversion formula then returns
$\arctan(x)$ at each $x\in\R$.
\end{example}

\begin{theorem}\label{theorempolynomial}
Let $f\in\bv([a,b])$ for each $-\infty<a<b<\infty$.  Suppose there are polynomials
$p_{\pm}(x)=\sum_{k=0}^{n_{\pm}}a_{\pm k} x^k$ such that
$f(x)\sim p_{\pm}(x)$ as $x\to\pm\infty$.  And, if
$
g(x)=f(x)-H(x)p_+(x)-H(-x)p_-(x)
$
then $g\in\bv(\R)$ with $\lim_{\abs{x}\to\infty}g(x)=0$.  Then, for $s\not=0$,
$$
\fhat(s)={\hat g}(s)+\sum_{k=0}^{n_+}a_{+k}\left[\pi i^k\delta^{(k)}(s)+\frac{k!}{(is)^{k+1}}\right]
+\sum_{k=0}^{n_-}a_{-k}\left[\pi i^k\delta^{(k)}(s)-\frac{k!}{(is)^{k+1}}\right].
$$
And, for each $x\in\R$,
\begin{equation}
\frac{f(x-)+f(x+)}{2}=\frac{1}{2\pi}\lim_{\substack{T\to\infty\\S\to 0^+}}\int_{S<\abs{s}<T} 
e^{ixs}\hat{g}(s)\,ds+H(x)p_+(x)+H(-x)p_-(x).
\label{polynomialinversion}
\end{equation}
\end{theorem}

\begin{proof}
Let $h_n(x)=H(x)x^n$.  The distributional Fourier transform

$$
{\hat h_n}(s)=\pi i^n\delta^{(n)}(s) +\frac{n!}{(is)^{n+1}}
$$
is well-known.  For example, \cite[p.~172]{gelfandshilov}. 
Let $q(x)=p_-(-x)$.  Then
$\fhat(s)={\hat g}(s)+{\widehat{Hp_+}}(s)+ {\widehat{Hq}}(-s)$.
And, 
$$
{\widehat{Hp_+}}(s)  =  \sum_{k=0}^{n_+}a_{+k}{\widehat h_n}(s)=
\sum_{k=0}^{n_+}a_{+k}\left[\pi i^k\delta^{(k)}(s)+\frac{k!}{(is)^{k+1}}\right].
$$
The transform of $Hq$ is similar, using the fact that $\delta^{(k)}(-s)=(-1)^k\delta^{(k)}(s)$.
The rest of the proof follows as with Theorem~\ref{theorembvlimit}.
\end{proof}

\begin{example}\label{examplepoly}
Let $p$ be a polynomial. Let $f(x)=p(x)\tanh(x)$ where $\tanh(x)=\sinh(x)/\cosh(x)$.  
Define 
$
g(x)= f(x)-{\rm sgn}(x)p(x)$.
Then as $\abs{x}\to\infty$ we have $g(x)\sim -2\,{\rm sgn}(x)p(x)e^{-2\abs{x}}$.  Hence, $g\in\Lone$.  
It follows that ${\hat g}(s)$ is continuous
on $\R$ and has limit $0$ as $\abs{s}\to\infty$.
Formula \eqref{polynomialinversion} gives
$f(x)$ for each $x$.  If $p(x)=O(x^2)$ as $x\to 0$ then $g$ and $g'$ are absolutely continuous and 
$g,g',g''\in\Lone$.  Integration by parts then shows
${\hat g}\in\Lone$.  Then for each $x\in\R$ the 
inversion formula reads
$$
f(x)=\frac{1}{2\pi}\int_{-\infty}^\infty e^{ixs}{\hat g}(s)\,ds+{\rm sgn}(x)p(x).
$$
\end{example}

\section{Principal value transforms}\label{sectionpv}

The inversion formula \eqref{localinversion} can be combined with conditions on $f$ outside
the interval $[a,b]$.  If $[b,\infty)$ is written as a finite union of intervals
$[b_0,b_1]\cup[b_1,b_2]\cup\ldots\cup[b_{N-1},b_N]$ (with $b=b_0<b_1<\ldots< b_N=\infty$) then to complete an 
inversion formula such as \eqref{inversionformula} we require
$$
\lim_{T\to\infty}\int_{-T}^T e^{ixs}\int_{b_{i-1}}^{b_i}e^{-ist}f(t)\,dt\,ds=0
$$
for each $1\leq i\leq N$.  Similarly on $(-\infty,a]$.  Use of the Dirichlet kernel in the proof of 
Corollary~\ref{corollary1}
shows $f\in L^1((b_{i-1},{b_i}))$ 
suffices.  The following theorem shows that on finite intervals we can replace this Lebesgue integration
condition with a principal value condition.

\begin{theorem}\label{theorempv}
Suppose $f$ is odd about $c$ such that $\int_c^{c+\delta}\abs{f(t)}\abs{t-c}\,dt$
exists for some $\delta>0$.  Define
the Fourier transform of $f$ by the principal value integral
$\fhat(s)=\int_{c-\delta}^{c+\delta}e^{-ist}f(t)\,dt$.
Suppose $x\not\in[c-\delta,c+\delta]$ and define $T_n=(2n+1)\pi/(2\abs{x-c})$ for $n\in\N$.
Then 
\begin{equation}\label{pv}
\lim_{\substack{n\to\infty\\S\to 0^+}}\int_{S<\abs{s}<T_n} 
e^{ixs}\fhat(s)\,ds =0.
\end{equation}
\end{theorem}

Note that $T_n\to\infty$ only on a special sequence depending on $x$ while $S\to 0^+$ unrestrictedly
through
real numbers.  M.~Mikol\'as \cite{mikolas} has also considered limits through a special sequence.
N.~Ortner has defined principal value Fourier transforms as distributions and also considers slowly
growing distributions \cite{ortner}.

\begin{proof}
By translation and reflection we see there is no loss of generality in assuming $f$ is
odd about $0$ and the integral $\int_0^\delta \abs{tf(t)}\,dt$ exists.  The
principal value Fourier transform is then
$\fhat(s)=-2i\int_0^\delta\sin(st)f(t)\,dt$.  Suppose $\abs{x}>\delta$.   Then use the
Fubini--Tonelli theorem and a trigonometric identity to write
\begin{align}
\int_{S<\abs{s}<T_n} 
&e^{ixs}\fhat(s)\,ds  =-2i\left\{\int_{-T_n}^{-S}+\int_{S}^{T_n}\right\}e^{ixs}\int_0^\delta
\sin(st)f(t)\,dt\,ds\notag\\
&=  4\int_S^{T_n} \sin(xs)\int_0^\delta\sin(st)f(t)\,dt\,ds\notag\\
& =  4\int_0^\delta f(t)\int_S^{T_n} \sin(xs)\sin(st)\,ds\,dt\notag\\
& =  2\int_0^\delta tf(t)\left[\frac{\phi_{T_n}(t-x)-\phi_S(t-x)-\phi_{T_n}(t+x)+\phi_S(t+x)}{t}\right]\,dt,
\label{deltaintegral}
\end{align}
where $\phi_U(u)=\sin(uU)/u$.  

We have
$$
\frac{\phi_{T_n}(t-x)-\phi_{T_n}(t+x)}{t}  =  \frac{2[t\sin(xT_n)\cos(tT_n)-x\cos(xT_n)\sin(tT_n)]}{(x^2-t^2)t}.
$$
Note that $\sin(tT_n)/t\sim T_n$ as $t\to 0$.
Since this is not bounded as a function of
$T_n$, we cannot, in general, bring the limit $T_n\to\infty$ inside the integral if $T_n$ is a
real variable.  But, with the choice $T_n=(2n+1)\pi/(2\abs{x})$
we have
$$
\frac{\phi_{T_n}(t-x)-\phi_{T_n}(t+x)}{t}  =  \frac{2\sin(xT_n)\cos(tT_n)}{x^2-t^2}
=\frac{2\,{\rm sgn}(x)(-1)^n}{x^2-t^2}\cos\left(\tfrac{(2n+1)\pi t}{2x}\right).
$$
Then, by the Riemann--Lebesgue lemma,
\begin{align*}
&\lim_{n\to\infty }\int_0^\delta tf(t)\left[\frac{\phi_{T_n}(t-x)-\phi_{T_n}(t+x)}{t}\right]\,dt\\
&  =   
2\,{\rm sgn}(x)\lim_{n\to\infty }(-1)^n\int_0^\delta\frac{tf(t)}{x^2-t^2}\cos\left(\tfrac{(2n+1)\pi t}{2x}\right)
\,dt=0.
\end{align*}

Similarly,
\begin{eqnarray*}
\left|\frac{-\phi_{S}(t-x)+\phi_{S}(t+x)}{t}\right| & = & 
\left|\frac{2[x\sin(tS)\cos(xS)-t\cos(tS)\sin(xS)]}{(x^2-t^2)t}\right|\\
 & \leq & \frac{2(\abs{x}S+1)}{x^2-\delta^2}.
\end{eqnarray*}
By dominated convergence,
$$
\lim_{S\to 0^+}\int_0^\delta tf(t)\left[\frac{-\phi_S(t-x)+\phi_S(t+x)}{t}\right]\,dt=0
$$
since
we can take the limit $S\to 0^+$ inside the integral and
$\phi_S(t\pm x)\to 0$ as $S\to 0$.
\end{proof}

\begin{example}
Let $f(x)={\rm sgn}(x)\abs{x}^{-\alpha}$ for $0<\alpha<2$.  Then the principal value
Fourier transform is
$\fhat(s)=-2i\,{\rm sgn}(s)\abs{s}^{\alpha-1}\int_0^\infty \sin(t) t^{-\alpha}\,dt$.
Note that $\fhat$ is Lebesgue integrable at the origin but integration in a neighbourhood
of infinity requires the inversion be 
a principal value integral.  
If $x\not=0$ we can
combine Corollary~\ref{corollarypringsheim} and Theorem~\ref{theorempv} to get the pointwise
inversion
$$
f(x)  = 
\frac{1}{2\pi}\lim_{n\to\infty}\int_{-T_n}^{T_n}
e^{ixs}\fhat(s)\,ds,
$$
where $T_n=(2n+1)\pi/(2\abs{x})$.

Take $0<\epsilon<\abs{x}$.
Write $f=f_1+f_2$ where $f_1=f\chi_{(-\epsilon,\epsilon)}$.  Then
\begin{eqnarray*}
{\hat f_1}(s) & = & \frac{2}{i}\int_0^\epsilon\frac{\sin(st)}{t^\alpha}\,dt
=\frac{2}{i}{\rm sgn}(s)\abs{s}^{\alpha-1}\int_0^{\epsilon\abs{s}}\frac{\sin(t)}{t^\alpha}\,dt\\
 & \sim & \frac{2}{i}\frac{\epsilon^{2-\alpha} s}{2-\alpha} \quad\text{ as } s\to0\\
{\hat f_1}(s) & \sim & \frac{2}{i}{\rm sgn}(s)\abs{s}^{\alpha-1}\int_0^\infty \frac{\sin(t)}{t^\alpha}\,dt 
\quad\text{ as } \abs{s}\to\infty.\\
\end{eqnarray*}
Hence, ${\hat f_1}$ is continuous at the origin but the principal value inversion of Theorem~\ref{theorempv}
is required since ${\hat f_1}$ is not integrable in a neighbourhood of infinity. 
Integration by parts and the Riemann--Lebesgue lemma show
$$
\lim_{n\to\infty}\int_{-T_n}^{T_n}
e^{ixs}{\hat f_1}(s)\,ds =0.
$$

And, 
\begin{eqnarray*}
{\hat f_2}(s) & = & \frac{2}{i}\int_\epsilon^\infty\frac{\sin(st)}{t^\alpha}\,dt
=\frac{2}{i}{\rm sgn}(s)\abs{s}^{\alpha-1}\int_{\epsilon\abs{s}}^\infty\frac{\sin(t)}{t^\alpha}\,dt\\
 & \sim & \frac{2}{i}{\rm sgn}(s)\abs{s}^{\alpha-1}\int_0^\infty \frac{\sin(t)}{t^\alpha}\,dt \quad\text{ as } s\to0\\
{\hat f_2}(s) & \sim & \frac{2}{i}\frac{\cos(\epsilon s)}{\epsilon^\alpha s} 
\quad\text{ as } \abs{s}\to\infty.
\end{eqnarray*}
In this case, ${\hat f_2}$ is Lebesgue integrable at the origin but the inversion integral will exist
conditionally.  Corollary~\ref{corollarypringsheim} gives
$$
f_2(x)=\frac{1}{2\pi}\lim_{T\to\infty}\int_{-T}^T
e^{ixs}{\hat f_2}(s)\,ds.
$$

Let $g(x)={\rm sgn}(x)\sin(ax)\abs{x}^{-\alpha}$ for $a\in\R$.  Then ${\hat g}(s)=-i(\fhat(s-a)-\fhat(s+a))/2$.
Let $h(x)={\rm sgn}(x)\cos(ax)\abs{x}^{-\alpha}$.   Then ${\hat h}(s)=(\fhat(s-a)+\fhat(s+a))/2$.
There are similar inversion formulas.
\end{example}

\section{The Henstock--Stieltjes integral}\label{sectionHSintegral}

This is an integral that properly extends the Riemann--Stieltjes integral.  It is described in 
detail on
the real line in
\cite{mcleod} and on compact intervals in \cite{monteiroslaviktvrdy} and \cite{tvrdy}.
The name Henstock is sometimes replaced with Kurzweil, Denjoy
or Perron.  Let $-\infty\leq a<b\leq\infty$.  A tagged partition
is a finite collection ${\mathcal P}=\{([x_{n-1},x_n],z_n)\}_{n=1}^N$ with tag $z_n\in[x_{n-1},x_n]$ and
$a=x_0<x_1<\ldots<x_N=b$, for some $N\in\N$.  A gauge is a function $\gamma$ from $[a,b]$ to the open
intervals in $[a,b]$ where we consider intervals such as $[a,x)$, $(x,y)$ and
$(x,b]$ to be open for $x,y\in(a,b)$.  (For $a=-\infty$ and $b=\infty$ this provides a two-point
compactification of the real line.)  Then ${\mathcal P}$ is said to be $\gamma$-fine if for each $1\leq n\leq N$
we have $[x_{n-1},x_n]\subset\gamma(z_n)$.  Let $f,g\fn[a,b]\to\R$.  Then $f$ is
integrable with respect to $g$ with integral $\int_a^b f\,dg\in\R$ if for each $\epsilon>0$ 
there is a gauge $\gamma$ such that if the tagged partition ${\mathcal P}$
is $\gamma$-fine then
\begin{equation}
\left|\sumN f(z_n)[g(x_n)-g(x_{n-1})]-\int_a^b f\,dg\right|<\epsilon.\label{HSdefn}
\end{equation}
For every $\epsilon>0$ it is always possible to choose a gauge $\gamma$ so that each $\gamma$-fine
tagged partition must have $a$ and $b$ as tags (of exactly one interval each).  We will always assume such a gauge has been chosen.

The Henstock--Stieltjes integral differs from the Riemann--Stieltjes integral in that it can integrate over
unbounded intervals and the tags, $z_n$, are chosen more carefully.  This allows the functions $f$ and $g$ to have
a common point of discontinuity, in which case the Riemann--Stieltjes integral would not
exist.

\begin{example}\label{exampleRS}
We verify
\eqref{regulatedidentity} for Henstock--Stieltjes integrals and show the integral can fail
to exist as a Riemann--Stieltjes integral, even allowing for improper integrals.

Define 
$$
f(x)=\left\{\begin{array}{cl}
0, & x<0\\
a, & x=0\\
1, & x>0,
\end{array}
\right.
$$
for some $a\in\R$.  Take $x=0$ in \eqref{regulatedidentity}.  The left side of this
identity is then $(f(0+)+f(0-))/2=1/2$.

Note that $H(-t)=0$
for $t>0$, $df=0$ except in a neighbourhood of the origin and $f=0$ on $(-\infty,0)$.

From \eqref{regulatedidentity}, $-i\omega\intinf H(x-t) e^{i\omega(x- t)}f(t)\,dt$
reduces to $-i\omega\int_{-\infty}^0 e^{-i\omega t}f(t)\,dt=0$, the integral existing
in the Lebesgue sense.

The only tag that can contribute to the Riemann sum in \eqref{HSdefn} is $z_n=0$.
Let $\epsilon>0$. Let $\delta>0$.  Define a gauge by $\gamma(0)=(-\delta,\delta)$.
For each $x\neq 0$ let $\gamma(x)$ be an open interval that contains $x$ but does not contain
the origin.  Define $\gamma(-\infty)=[-\infty,-M)$ and $\gamma(\infty)=(M,\infty]$ for
some $M>0$.  Combining tagged intervals with coincident endpoints if necessary,
every $\gamma$-fine tagged partition of $[-\infty,\infty]$ must then
have exactly one tagged interval with $0$ as a tag.  This is the only tag that
can contribute to the Riemann sum representing $\intinf H(x-t) e^{i\omega(x- t)}\,df(t)$.
Write this tagged interval as $([\alpha,\beta],0)$ where $-\delta<\alpha<0<\beta<\delta$.
The Riemann sum in \eqref{HSdefn} then reduces to $H(0)[f(\beta)-f(\alpha)]=1/2$.  This
verifies \eqref{regulatedidentity} for Henstock--Stieltjes integrals.

Now consider $\intinf H(-t) e^{-i\omega t}\,df(t)$ as a Riemann--Stieltjes integral.
As noted above, the integrand or integrator vanish except in a neighbourhood of the origin.
We can then just integrate over $[-1,1]$.  With Riemann--Stieltjes integrals, the length
of intervals allowed in a Riemann sum are controlled by a parameter.  Suppose $\delta>0$
and all intervals
must have length less than $\delta$.  In a partition of $[-1,1]$ we 
are free to choose an interval $[\alpha,\beta]$
where $-1\leq\alpha<0<\beta\leq 1$ and $\beta-\alpha<\delta$.  We are free to choose any tag
$z\in[\alpha,\beta]$.  As above, this is the only tag that contributes non-trivially to
a Riemann sum.  The Riemann sum then reduces to $H(z)e^{-i\omega z}[f(\beta)-f(\alpha)]=
H(z)e^{-i\omega z}/2$.
If $z=\alpha$ this gives $0$.  If $z=\beta$ this gives $e^{-i\omega \beta}/2$.  Hence, the Riemann--Stieltjes
integral does not exist.

The improper Riemann--Stieltjes integral also does not exist.  If $0<\epsilon<1$ then
$\int_\epsilon^1 H(-t) e^{-i\omega t}\,df(t)$ and 
$\int_{-1}^{-\epsilon} H(-t) e^{-i\omega t}\,df(t)$ both vanish since $f$ is constant
on the intervals of integration.  Hence, any improper integral computed by taking
the limit $\epsilon\to 0^+$ will yield $0$, not the required value of $1/2$.
\end{example}

If one of $f$ and $g$ is regulated and the other is of bounded variation then the integrals $\int_a^bf\,dg$
and $\int_a^bg\,df$ exist and they are related by an integration by parts formula.  
It allows the functions to have discontinuities at the same point.
The following formula is proved in \cite[p.~199]{mcleod}.  Let $-\infty\leq a< b\leq\infty$.  Then

\begin{align}
\int_a^bf\,dg  =&  f(b)g(b)-f(a)g(a)-\int_a^b g\,df
+ \sum[f(c_n)-f(c_n-)][g(c_n)-g(c_n-)]\notag\\
&\qquad-\sum[f(c_n)-f(c_n+)][g(c_n)-g(c_n+)].
\label{parts}
\end{align}
The sums contain all points $c_n$ at which $f$ and $g$ are either discontinuous from the left or 
discontinuous from 
the right.

We now prove two related results.
\begin{prop}\label{prop1}
Let $-\infty\leq a<b\leq\infty$.  Let $A,B,C\fn[a,b]\to\R$ such that two functions are in $\bv([a,b])$ and the
other is in $\reg([a,b])$.  If $B$ and $C$ have no common point of discontinuity then
$$
\intab A\,d[BC] =\intab AB\,dC+\intab AC\,dB.
$$
\end{prop}
\begin{proof}
Since $\bv([a,b])$ and $\reg([a,b])$ are both closed under pointwise products, the integrals in the proposition
exist.

Given $\epsilon>0$ there is a gauge $\gamma$ and a
$\gamma$-fine tagged partition of $[a,b]$, defined by $\{[x_{n-1},x_n],z_n\}_{n=1}^N$, such that
\begin{align*}
&\left|\sumN A(z_n)[B(x_n)C(x_n)-B(x_{n-1})C(x_{n-1})]-\intab A\,d[BC]\right|<\epsilon\\
&\left|\sumN A(z_n)B(z_n)[C(x_n)-C(x_{n-1})]-\intab AB\,dC\right|<\epsilon\\
&\left|\sumN A(z_n)C(z_n)[B(x_n)-B(x_{n-1})]-\intab AC\,dB\right|<\epsilon\\
&\left|\sumN B(z_n)[C(x_n)-C(x_{n-1})]-\intab B\,dC\right|<\epsilon\\
&\left|\sumN C(z_n)[B(x_n)-B(x_{n-1})]-\intab C\,dB\right|<\epsilon\\
&\left|\sumN [B(x_n)C(x_n)-B(x_{n-1})C(x_{n-1})]-\intab d[BC]\right|<\epsilon.
\end{align*}
By Henstock's lemma \cite[p.~186]{mcleod},
\begin{align*}
&\sumN\left|B(z_n)[C(x_n)-C(x_{n-1})]-\int_{x_{n-1}}^{x_n} B\,dC\right|<2\epsilon\\
&\sumN\left|C(z_n)[B(x_n)-B(x_{n-1})]-\int_{x_{n-1}}^{x_n} C\,dB\right|<2\epsilon.
\end{align*}
Let
\begin{align*}
D_n=&B(x_n)C(x_n)-B(x_{n-1})C(x_{n-1})-B(z_n)[C(x_n)-C(x_{n-1})]\\
&\qquad-C(z_n)[B(x_n)-B(x_{n-1})].
\end{align*}
It then suffices to show $\sum_{n=1}^N \abs{A(z_n)D_n}$ is less than a multiple of $\epsilon$.  Notice that 
since $B$ and $C$ have no common points of discontinuity, the integration by parts formula is
$$
\int_{x_{n-1}}^{x_n} d[BC]=B(x_n)C(x_n)-B(x_{n-1})C(x_{n-1})=\int_{x_{n-1}}^{x_n} B\,dC+\int_{x_{n-1}}^{x_n} C\,dB.
$$
The above inequalities then give
$$
\sum_{n=1}^N \abs{D_n}  =  \sumN\left|D_n-\int_{x_{n-1}}^{x_n} d[BC]+\int_{x_{n-1}}^{x_n} B\,dC+\int_{x_{n-1}}^{x_n} C\,dB\right|
  <  4\epsilon.
$$
Then $\sum_{n=1}^N \abs{A(z_n)D_n}<4\epsilon\norm{A}_\infty$.
\end{proof}

If $B$ and $C$ have coincident discontinuities at $c_n$ then sum terms as in the integration by
parts formula need to be added: $\sum A(c_n)[B(c_n)-B(c_n-)][C(c_n)-C(c_n-)]$ and
$\sum A(c_n)[B(c_n)-B(c_n+)][C(c_n)-C(c_n+)]$.

\begin{prop}\label{propB'}
Let $A\in\reg([a,b])$ and $B$ be absolutely continuous such that $B'\in L^1([a,b])$.  Then
$\intab A\,dB=\intab A(t)B'(t)\,dt$.
\end{prop}
\begin{proof}
If $[a,b]$ is a finite interval then $B\in\bv([a,b])$.  If $[a,b]$ is not finite we can 
write
$$
\sumN\abs{B(x_n)-B(x_{n-1})}=\sumN\left|\int_{x_{n-1}}^{x_n}B'(t)\,dt\right|\leq
\sumN\int_{x_{n-1}}^{x_n}\abs{B'(t)}\,dt=\norm{B'}_1.
$$
Hence, $B\in\bv([a,b])$. 
The rest of the proof is similar to that for Proposition~\ref{prop1}, starting with a gauge that
simultaneously makes Riemann sums within $\epsilon$ of the integrals
$\intab A\,dB$, $\intab A(t)B'(t)\,dt$ and $\intab B'(t)\,dt$.
\end{proof}
Notice that in the case of unbounded intervals absolute continuity of $B$ does not imply
$B'\in L^1([a,b])$.  For example, $B(x)=x$.

\section{Proof of Lemma~\ref{lemmahsrep}}\label{lemmaproofhsrep}

\begin{proof}
(a) Using the cutoff property of the Heaviside step function, the integration parts formula 
\eqref{parts} and the fact that $f$ is bounded gives
\begin{align*}
&\intinf H(x-t) e^{-i\omega t}[df(t)-i\omega f(t)\,dt]
=\int_{-\infty}^b H(x-t) e^{-i\omega t}[df(t)-i\omega f(t)\,dt]\\
&\quad=H(x-b)e^{-i\omega b}f(b)
-\lim_{t\to-\infty}\left[H(x-t)e^{-i\omega t}f(t)\right]
-\int_{-\infty}^b f(t)\,d\left[H(x-t)e^{-i\omega t}\right]\\
&\quad+e^{-i\omega x}\left\{\left[H(0)-H(0+)\right]\left[f(x)-f(x-)\right]
-\left[H(0)-H(0-)\right]\left[f(x)-f(x+)\right]\right\}\\
&\quad-i\omega\int_{-\infty}^b H(x-t)e^{-i\omega t}f(t)\,dt\\
&=-\int_{-\infty}^b f(t)e^{-i\omega t}\,dH(x-t)
-e^{-i\omega x}\left\{\left[\frac{f(x)-f(x-)}{2}\right]
+\left[\frac{f(x)-f(x+)}{2}\right]\right\}\\
&=e^{-i\omega x}\left[\frac{f(x-)+f(x+)}{2}\right].
\end{align*}

Since the function $t\mapsto e^{-i\omega t}$ is absolutely continuous on $(-\infty,b]$
and its derivative is in $L^1((-\infty,b])$
we can use
Propositions~\ref{prop1} and \ref{propB'} in the last lines above.
A gauge can be chosen so that each 
$\gamma$-fine tagged partition must have $x$ as a tag.  This shows
$\int_{-\infty}^b f(t)e^{-i\omega t}\,dH(x-t)=f(x)e^{-i\omega x}[H(0-)-H(0+)]$.

Part (b) is similar.
\end{proof}

\end{document}